\documentclass[11pt]{amsart}
\usepackage{amssymb}
\usepackage{amscd}
\theoremstyle{plain}
\newtheorem{theorem}{Theorem}
\theoremstyle{definition}
\newtheorem{proposition}{Proposition}[section]
\newtheorem{corollary}[proposition]{Corollary}
\newtheorem{lemma}[proposition]{Lemma}
\newtheorem{definition}[proposition]{Definition}

\theoremstyle{remark}

\newcommand{\sk}{\vskip.05in}

\newcommand{\restr}{\upharpoonright}

\newcommand{\subs}{\subseteq}

\numberwithin{equation}{section}
\begin{document}
\title[Separable monotonically normal compacta]{Elementary submodels and separable monotonically normal compacta}
\author{Todd Eisworth}
\address{Department of Mathematics\\
         Ohio University\\
         Athens, OH 45701}
\email{eisworth@math.ohiou.edu}
\date{\today}
\begin{abstract}
In this note, we use elementary submodels to prove that a separable
monotonically normal compactum can be mapped onto a separable metric
space via a continuous function whose fibers have cardinality at
most 2.
\end{abstract}
\thanks{The author acknowledges support from NSF grant DMS-0506063.}
\maketitle
\section{Introduction}
Mary Ellen Rudin's 2001 proof of Nikiel's Conjecture \cite{nikiel} is a tour de force of combinatorial
 set theory and general topology. The work in this note arose from the author's attempt at understanding
  her proof; it became clear upon reading her paper that elementary submodels ought to be able to shed light
  on the structure of monotonically normal compacta. In the sequel, we will give one such application.

We will begin with the definition of monotonically normal spaces.

\begin{definition}
A space\footnote{We assume that all topological spaces under consideration in the paper are Tychonoff.}
 $X$ is monotonically normal if there is a binary function $H$ whose domain consists of all pairs
 $(p, U)$ where $p\in U$ and $U$ is open in $X$ such that
 \begin{enumerate}
 \sk
 \item $H(p, U)$ is an open set such that $p\in H(p, U)\subs U$,
 \sk
 \item if $p\notin V$ and $q\notin U$ then $H(p, U)\cap H(p, V)=\emptyset$, and
 \sk
 \item if $V\subs U$ is another open set, then $H(p, U)\subs H(p, V)$.
 \end{enumerate}
 The function $H$ is called a {\em monotone normality operator} on $X$.
 \end{definition}

It is easy to see that a monotonically normal space $X$ is normal --- if $H$ and $K$ are disjoint closed
subsets of $X$, and we define
\begin{equation}
U=\bigcup_{p\in H} H(p, X\setminus K)
\end{equation}
and
\begin{equation}
V=\bigcup_{q\in K} H(q, X\setminus H),
\end{equation}
then $U$ and $V$ provide the required separation of $H$ and $K$.

There is an extensive literature devoted to monotonically normal spaces (Gruenhage's two articles
\cite{gruenhageold} and \cite{gruenhagenew} can be consulted for more information); much of
Rudin's contribution to this literature centered on Nikiel's Conjecture, and her work culminated
in the following theorem.

\begin{theorem}[Nikiel's Conjecture -- Rudin \cite{nikiel}]
A compact Hausdorff space $X$ is monotonically normal if and only if
it is the continuous image of a linearly ordered compactum.
\end{theorem}

Our work in this area begins with her proof, but before we can prove our theorem, we must deal with some
aspects of elementary submodels in topology.

\section{Elementary submodels and $X/M$}

We will assume that the reader has some familiarity with the use of elementary submodels in topology ---
Dow's papers \cite{dow88} and \cite{dow95} provide more than adequate preparation.  We are going to be using
a particular construction involving elementary submodels due independently to Bandlow (\cite{bandlow1},
 \cite{bandlow2}, and \cite{bandlow3}) and Dow \cite{dow92}, and this section of the paper will be used to develop the theory in a self--contained manner.

Let us assume now that $(X,\tau)$ is a Tychonoff space, $\chi$ is a ``sufficiently large'' regular cardinal,
and $M$ is an elementary submodel of $H(\chi)$, the collection of sets hereditarily of cardinality $<\chi$.

\begin{definition}
\label{defn1}
Two points $x$ and $y$ are $M$--equivalent, denoted $x=_M y$,  if $f(x)=f(y)$ for every continuous function $f:X\rightarrow\mathbb{R}$
such that $f\in M$.
\end{definition}

\begin{proposition}
$M$--equivalence is an equivalence relation, and the equivalence classes are closed subsets of $X$.
\end{proposition}
\begin{proof}
The fact that $=_M$ is an equivalence relation is trivial.  To see that equivalence classes are closed,
suppose that $x$ is not $M$--equivalent to $y$, and let $f\in M$ be a function that witnesses this.
Then $f^{-1}(\mathbb{R}\setminus\{f(y)\})$ is an open neighborhood of $x$ that is disjoint to the $M$--equivalence
class of $y$.
\end{proof}

We denote the equivalence class of $x$ by $[x]$, and let $X/M$ denote the set of all $M$--equivalence classes
of points in $X$. There are several natural choices for topologizing $X/M$ (all of which coincide if $X$ is compact),
but the theory works smoothest when then following definition is used.

\begin{definition}
Let $\pi_M$ be the natural projection of $X$ onto $X/M$.  We topology $X/M$ by taking as a base all sets of the form
$\pi_M[U]$ for $U$ a co-zero set in $M$.
\end{definition}

The proof of the following proposition is left to the reader --- all of the clauses enumerated follow
easily from the definitions involved.

\begin{proposition}
Let $X$ and $M$ be as above.
\begin{enumerate}
\item The mapping $\pi_M:X\rightarrow X/M$ is continuous.
\sk
\item $x=_M y$ if and only if for every co--zero $U\in M$,
\begin{equation}
x\in U\Longleftrightarrow y\in U.
\end{equation}
\item $[x]=\bigcap\{Z\in M: x\in Z\text{, $Z$ a zero--set}\}$.
\sk
\item $[x]=\bigcap\{\overline{U}: x\in U\in M\text{, $U$ a co--zero set}\}.$
\sk
\item $X/M$ is a Hausdorff space.
\sk
\end{enumerate}
\end{proposition}

We will prove shortly that $X/M$ is in fact a Tychonoff space, but before we do that we point out the
connection between $X/M$ and a more typical construction using elementary submodels.

\begin{definition}
Given $X$ and $M$ as above, $X_M$ is defined to be the topological space whose underlying set is $M\cap X$
with the topology generated by sets of the form $M\cap U$, where $U$ is an open subset of $X$ that is in $M$.
\end{definition}

\begin{proposition}
$X_M$ is homeomorphic to a dense subspace of $X/M$.
\end{proposition}
\begin{proof}
The proof is the obvious one --- a point $x\in M\cap X$ is set to its equivalence class in $X/M$.
The topologies involved are defined in such a way that this map is a homeomorphism of $X_M$ onto its
image, and the fact that this image is dense in $X/M$ follows easily as any non--empty open subset of $X$ that is
in $M$ must contain an element of $M\cap X$.
\end{proof}

Our next goal is to give a more concrete description of $X/M$ in terms of $\beta X$, the Stone-Cech compactification
of $X$.  One of the standard constructions of $\beta X$ involves embedding $X$ into a product $[0,1]^\kappa$ of
unit intervals, so we first investigate the nature of $X/M$ for $X$ of the form $[0,1]^\kappa$.

\begin{proposition}
Suppose $X= [0, 1]^\kappa$ for some cardinal $\kappa$ (so we view points of $X$ as functions from $\kappa$ to $[0,1]$),
 and let $M$ be an elementary submodel of $H(\chi)$  containing $X$. Then
\begin{enumerate}
\item $x =_M y$ if and only if $x\restr M\cap\kappa = y\restr M\cap\kappa$, and
\sk
\item $X/M$ is homeomorphic to $[0, 1]^{M\cap \kappa}$.
\sk
\end{enumerate}
\end{proposition}
\begin{proof}
It is certainly the case that $x\restr M\cap\kappa = y\restr M\cap\kappa$ if $x=_M y$, as projection onto
the $\alpha$th component is a real--valued function on $X$ that is in $M$ if $\alpha\in M\cap\kappa$. Suppose
now that $x\neq_M y$; we produce an $\alpha\in M\cap \kappa$ such that $x(\alpha)\neq y(\alpha)$. The key to
this is a well--known fact that a continuous real--valued function defined on a product of compact spaces
depends on countably many coordinates.\footnote{See Exercise 3.2H(a) of Engelking \cite{engelking}, for example.}
In our context, this means if we have a function $f:X\rightarrow\mathbb{R}$ then there is a countable
set $S\subs\kappa$ such that $f(x)=f(y)$ if $x\restr S = y\restr S$. If $x\neq_M y$, then there is a
function $f\in M$ mapping $X$ to $\mathbb{R}$ such that $f(x)\neq f(y)$. The model $M$ knows that $f$
depends on countably many coordinates, so there is a countable $S\in M$ with the requisite properties.
In particular, there is an $\alpha\in S$ such that $x(\alpha)\neq y(\alpha)$.  Since the set $S$ is countable
and an element of $M$, we know that $S\subs M$ and therefore there is an $\alpha\in S\subs M\cap\kappa$ for
which $x(\alpha)\neq y(\alpha)$, as required.

Thus, there is a natural correspondence between $X/M$ and $[0, 1]^{M\cap\kappa}$.  The fact that this correspondence
is a homeomorphism is not difficult to show, and is left to the reader.
\end{proof}

Now suppose we have $X$ and $M$ as usual. Since $\chi$ is ``large enough'', we know that the space
$C^*(X)$ of continuous functions from $X$ to $[0,1]$ is going to be an element of $M$, as it is definable
from $X$ using parameters available in $M$.  Thus, inside $M$ there will be an enumeration $\langle f_\alpha:\alpha<\kappa\rangle$
of $C^*(X)$.\footnote{There will be many such enumerations and it doesn't really matter which one we use, but
for definitiveness, we can consider $M$ to be an elementary submodel
of the expanded structure $\langle H(\chi), \in, <_\chi\rangle$ where $<_\chi$ is a fixed well--ordering
of $H(\chi)$. We can then use the well--ordering to pick the particular enumeration we use.}

We know that the function $e: X\hookrightarrow [0, 1]^\kappa$  that maps $x$ to
 $\langle f_\alpha(x):\alpha<\kappa\rangle$
embeds $X$ as a subspace of $[0,1]^\kappa$.  The following proposition establishes a similar connection
between $X/M$ and $[0, 1]^{M\cap\kappa}$.

\begin{proposition}
Suppose $X$, $M$, and $\langle f_\alpha:\alpha<\kappa\rangle$ are as in the preceding discussion.
Then there is an natural embedding $e/M$ making the following diagram commute:
\begin{equation}
\begin{CD}
X @>e>>  [0,1]^\kappa \\
@V\pi/MVV    @VV\pi/MV \\
X/M  @>>e/M> [0,1]^{M\cap\kappa}\\
\end{CD}
\end{equation}
\end{proposition}

The preceding material is taken from Section~5 of Dow's \cite{dow92}; the reader can find a more detailed
discussion there.  For our purposes, we need only the following corollary.

\begin{corollary}
$X/M$ is a Tychonoff space.
\end{corollary}
\begin{proof}
The embedding $e/M$ from the preceding proposition shows us that $X/M$ is homeomorphic to a subspace of
the compact space $[0, 1]^{M\cap\kappa}$.
\end{proof}

\section{Monotone Normality}

We are now in a position to state the main theorem of this note.

\begin{theorem}
\label{mainthm} Let $X$ be a separable monotonically normal
compactum, and let $M$ be a countable elementary submodel of
$H(\chi)$ containing $X$.  Then each $=_M$--equivalence class has
cardinality at most $2$.
\end{theorem}

As a corollary, we get the following result mentioned in the abstract.

\begin{corollary}
If $X$ is a separable monotonically normal compactum, then $X$ can
be continuously mapped onto a separable metric space by a 2-to-1
function\footnote{We operate under the convention that this means
pre-images of points have cardinality {\em at most} 2.}.
\end{corollary}
\begin{proof}
In light of Theorem~\ref{mainthm}, it suffices to note that $X/M$ is compact (as it is the continuous
image of $X$) and of countable weight (as the countably many co--zero sets in $M$ define a base).
\end{proof}

Actually, much more can be shown. In particular, $X$ can be written as the inverse limit of a system
of compact metric spaces in such a way that all of the projections from $X$ onto members of the system are
two-to-one functions.

Before we prove Theorem~\ref{mainthm}, we deal with the special case where $X$ is linearly ordered as this
is particularly simple, and it hints at why Theorem~\ref{mainthm} is connected to Nikiel's Conjecture.

Suppose now that $X$ is a separable linearly ordered space, and let $M$ be a countable elementary submodel
of $H(\chi)$ containing $X$.  In the model $M$, one can find a countable set $D\subs X$ that is dense in $X$.
Since $D$ is countable, it follows that every element of $D$ is also in $M$.

Now suppose $x<y$ are $M$--equivalent; it suffices to prove that the interval $(x, y)$ is empty.
Assume by way of contradiction that $(x, y)\neq\emptyset$. Then there is an element $d\in D$ such that
$x<d<y$. If one of $(x, d)$ and $(d, y)$ is empty, then we easily get a function in $M$ separating
$x$ and $y$ --- for example, if $(x, d)$ is empty then $x$ is in $M$ (it's definable in $M$ as the predecessor
of $d$), and hence so is the (continuous) function sending $(-\infty, x]$ to $0$ and $[d, \infty)$ to 1.
If both $(x, d)$ and $(d, y)$ are non--empty, then we can find $d_1$ and $d_2$ in $D$ such that
\begin{equation}
x<d_1<d_2<y
\end{equation}
and the disjoint closed sets $(-\infty, d_1]$ and $[d_2,\infty)$ are in $M$ and separate $x$ and $y$.  Since
$X$ is normal, we can find in $M$ a function separating these two closed sets, and the function also separates
$x$ and $y$.

Thus, if $X$ is separable and linearly ordered (we don't need compactness for this case) and $M$ is as above, then each $=_M$ equivalence class
is either a singleton, or a pair $x<y$ with $y$ the immediate successor of $x$.

\begin{proof}[Proof of Theorem~\ref{mainthm}]
Let $X$ be a separable monotonically normal compactum, and let $M$
be a countable elementary submodel of $H(\chi)$ that contains $X$.
Let $H\in M$ be a monotone normality operator.

\begin{lemma}
\label{lemma1}
If $K$ is an $=_M$--equivalence class and $U$ is any open neighborhood of $K$, then there is a co--zero
set $V\in M$ such that $K\subs V$ and $\overline{V}\subs U$.
\end{lemma}
\begin{proof}
We know that $K$ is the intersection of all co--zero sets in $M$ that contain it.  Since $X$ is compact
and $K$ is closed, this implies that the collection of co--zero sets from $M$ that contain $K$ is a base
for $K$, i.e., any open set containing $K$ must contain such a co--zero set.  The result follows immediately.
\end{proof}

\begin{lemma}
\label{lemma2}
If $K$ is an $=_M$--equivalence class and $U$ is any open neighborhood of $K$, then there is a point
$p\in M\cap U$ and a co--zero set $V\in M$ such that
\begin{enumerate}
\sk
\item $p\in M\cap V$,
\sk
\item $\overline{V}\subs U$, and
\sk
\item $K\subs H(p, V)$.
\end{enumerate}
\end{lemma}
\begin{proof}
As in the proof of Lemma~\ref{lemma1}, there is a co--zero set $V\in M$ such that $K\subs V$ and
$\overline{V}\subs U$.  Since $V$ is (in $M$) a countable union of zero sets, there is a zero set
$Z\in M$ such that $Z\cap K\neq\emptyset$ and therefore $K\subs Z$ by the definition of $M$--equivalence.
For each $p\in Z$ there is a co--zero set $V_p$ such that
\begin{equation}
p\in V_p\subs H(p, V).
\end{equation}
By elementarity, we may assume that the mapping $p\mapsto V_p$ is an element of $M$, and thus
$\{V_p:p\in Z\}$ is an open cover of $Z$ that is an element of $M$.  Since $Z$ is compact, there is
finite $Z_0\subs Z$ such that
\begin{equation}
Z\subs \bigcup_{z\in Z_0} V_p.
\end{equation}
We may assume that $Z_0\in M$, and therefore $Z_0\subs M$.  Thus, there is a $p\in M\cap Z$ such that
$V_p\cap K\neq\emptyset$.  Since $V_p$ is a co--zero set in $M$, it follows that $K\subs V_p$. Thus
\begin{equation}
K\subs V_p\subs H(p, V)\subs V,
\end{equation}
and the result follows.
\end{proof}

We now are in a position to apply the monotone normality of $X$ in a non--trivial way.

\begin{proposition}
\label{prop34}
Suppose $K$ is a $=_M$--equivalence class, $\{x_0, x_1\}\subs K$, and $W_0$ and $W_1$ are disjoint
open sets with $x_i\in W_i$.\footnote{So $W_0$ and $W_1$ won't be elements of $M$.}  If $K'$
is any other equivalence class, then there is at most one $i<2$ with $K\cap H(x_i, W_i)\neq\emptyset$.
\end{proposition}
\begin{proof}
Let $U$ and $U'$ be disjoint co--zero sets in $M$ separating $K$ and $K'$.  By the previous lemma,
there is a point $p\in M\cap U'$ such that
\begin{equation}
K'\subs H(p, U').
\end{equation}
Note that $p$ is an element of $W_i$ for at most one $i$, and that $K\cap U'=\emptyset$.  If
$p\notin W_i$, then
\begin{equation}
H(p, U')\cap H(x_i, W_i)=\emptyset
\end{equation}
because $H$ is a monotone normality operator.  Since $K'\subs H(p, U')$, it follows that
\begin{equation}
K'\cap H(x_i, W_i)=\emptyset
\end{equation}
as well.
\end{proof}

\begin{definition}
Let $U\in M$ be a co--zero set. An equivalence class $K\in X/M$ is said to be {\em shattered by $U$} if
$K\subs U$ and there exist $\{x_i:i<3\}$ and $\{W_i:i<3\}$ such that
\begin{itemize}
\sk
\item $\{x_i:i<3\}\subs K$,
\sk
\item $W_i$ is an open neighborhood of $x_i$,
\sk
\item the $W_i$'s have pairwise disjoint closures, and
\sk
\item $H^4(x_i, W_i)\nsubseteq U$ for all $i<3$.
\end{itemize}
(Here the notation $H^n(p, U)$ is defined by induction: $H^2(p, U)= H(H(p, U))$ and $H^{n+1}(p, U)=H(H^n(p, U))$.)
\end{definition}

We will show that each co--zero $U\in M$ can shatter at most countably many elements of $X/M$, and that
every $K\in X/M$ of cardinality greater than $2$ is shattered by some $U\in M$.

\begin{proposition}
\label{prop1}
Let $U\in M$ be a co--zero set. Then there are at most countably many equivalence classes in $X/M$ that
are shattered by $U$.
\end{proposition}
\begin{proof}
By way of contradiction, suppose there are uncountably many
equivalence classes shattered by $U$.  Since $U$ is a co--zero set
in $M$, we know there is a family $\{Z_n:n\in\omega\}$ of
zero--sets in $M$ whose union is $U$.  Recall that since $Z_n\in
M$, if $Z_n$ meets an equivalence class in $X/M$, then $Z_n$
actually contains the entire equivalence class.  Thus we can find
a zero--set $Z\in M$ and a family $\{K_n:n\in\omega\}$ of
equivalence classes in $X/M$ such that
\begin{itemize}
\item $Z\subs U$
\sk
\item $K_n\subs Z$ for all $n$
\sk
\item each $K_n$ is shattered by $U$
\end{itemize}
For $n<\omega$, let $\{x_i^n:i<3\}$ and $\{W^n_i:i<3\}$ be as in
the  previous definition for the equivalence class $K_n$.

For each pair of natural numbers $m<n$, by Proposition \ref{prop34} we
can find a value $i=i(m,n)<3$ such that
\begin{equation}
\label{eqn1}
 K_n\cap H(x^m_i, W^m_i)=\emptyset\text{ and } K_m\cap
H(x^n_i, W^n_i)=\emptyset.
\end{equation}
Therefore, by an application of Ramsey's Theorem, we may assume
that there is an $i$ such that (\ref{eqn1}) holds for all $m\neq
n$.

Since $H$ is a monotone normality operator, this implies that for
$m\neq n$,
\begin{equation}
H^2(x^n_i, W^n_i)\cap H^2(x^m_i, W^m_i)=\emptyset.
\end{equation}
For each $n$, choose $p_n\in H^4(x^n_i, W^N_i)\setminus U$.  Since
$X$ is compact, we can find a point $p$ that is a limit point of
$\{x_n:n\in\omega\}$.

Since the family $\{H^2(x^n_i, W^n_i):n<\omega\}$ is pairwise
disjoint, the members of $\{H^3(x^n_i, W^n_i:n<\omega\}$ have
pairwise disjoint closures.  Thus $p$ is a member of $H^3(x_i^n,
W^n_i)$ for at most one $n$.

If $p\notin H^3(x^m_i, W^m_i)$ then  $H(p, X\setminus U)$ and
$H^4(x^m_i, W^m_i)$ are disjoint, and so $H(p, X\setminus U)$ is
an open neighborhood of $p$ that contains at most one member of
$\{p_n:n<\omega\}$, a contradiction.
\end{proof}

\begin{proposition}
If $K\in X/M$ is an equivalence class of size $\geq 3$, then there
is a co--zero set $U\in M$ that shatters $K$.
\end{proposition}
\begin{proof}
Note that since $X$ is separable, an equivalence class in $X/M$
has empty interior, except in the case where the equivalence class
consists of a single isolated point from $X$. If $K\in X/M$ has
size $\geq 3$, choose distinct $\{x_i:i<3\}$ in $K$, and choose
open sets (not necessarily from $M$!) $\{W_i:i<3\}$ with disjoint
closures such that $x_i\in W_i$.

Since $K$ is the intersection of all co--zero sets from $M$ that
contain it and $K$ has empty interior, for each $i<3$ we can find
a cozero set $U_i\in M$ such that $K\subs U_i$ and $H^4(x_i,
W_i)\nsubseteq U_i$.  Finally, the set $U=U_0\cap U_1\cap U_2$ is
a cozero set in $M$ with all the required properties.

\end{proof}

From the two preceding propositions, it follows that all but countably many equivalence classes
 in $X/M$ are of size $2$ or smaller, but we need to improve this to {\em all} equivalence classes.  The proof
  of this breaks into two steps --- first we use a variant of an argument from Rudin's paper \cite{nikiel} to
  show that there is at least one $M$ for which every $=_M$--class has size at most 2, and then we
   show that in fact it must hold for every such $M$.

\begin{proposition}
\label{prop2}
There is a countable elementary submodel $M$ of $H(\chi)$ containing $X$ such that all $=_M$--classes are
of cardinality $\leq 2$.
\end{proposition}
\begin{proof}
Assume by way of contradiction that the result fails, and let $\langle M_\alpha:\alpha<\omega_1\rangle$
be an increasing and continuous $\in$--chain of countable elementary submodels of $H(\chi)$ such that
$X\in M_0$ and $\langle M_\beta:\beta\leq\alpha\rangle\in M_{\alpha+1}$.

By our assumption, we can choose for each $\alpha<\omega_1$ an $=_\alpha$--class\footnote{We write $=_\alpha$ instead of $=_{M_\alpha}$.} $K_\alpha\in M_{\alpha+1}$
containing at least three elements.  Since $M_\alpha\in M_\beta$ for $\alpha<\beta$, it follows that every
$=_\beta$--class is contained in a unique $=_\alpha$--class.  Thus,
\begin{equation}
\alpha<\beta<\omega_1\Longrightarrow\text{ either }K_\alpha\cap K_\beta=\emptyset\text{ or }K_\beta\subs K_\alpha.
\end{equation}
Since $\omega_1\rightarrow (\omega_1,\omega)^2$ by the Dushnik--Miller Theorem\footnote{From \cite{DM}, or see Theorem~14.6 of \cite{hh}}
we know that  either there is an infinite $A\subs\omega_1$ such that $K_\alpha\cap K_\beta=\emptyset$ for
$\alpha<\beta$ in $A$, or there is an uncountable set $B\subs\omega_1$ such that $K_\beta\subs K_\alpha$ for
$\alpha<\beta$ in $B$.  We will get a contradiction by proving that both of these alternatives are untenable.

We will dispose of the second alternative first, so assume that we have such an uncountable $B$.  Are argument
consists of showing that $K_\beta$ is actually a {\em proper} subset of $K_\alpha$ for $\alpha<\beta<\omega_1$, and
then quoting an old result due both to Ostaszewski and to Moody.

Suppose now that $\alpha<\beta$ and $K_\beta\subs K_\alpha$.  Since $K_\alpha\in M_{\alpha+1}$ and $M_{\alpha+1}$
is countable, we know that $K_\alpha\in M_\beta$ as well, and hence
\begin{equation}
M_\beta\models |K_\alpha|\geq 3.
\end{equation}
In particular, we can find points $x\neq y$ in the set $M_\beta\cap K_\alpha$. These two points are separated
by a continuous function in $M_\beta$, and hence $K_\beta$ can contain at most one of them. Thus, $K_\beta$
is a proper subset of $K_\alpha$.

This implies that the sequence $\langle K_\alpha:\alpha\in B\rangle$ is an uncountable strictly decreasing
sequence of closed subsets of $X$. However, this is absurd, as by Ostaszewski \cite{osta} and Moody \cite{moody}
a separable monotonically normal space is hereditarily Lindel\"{o}f.\footnote{In fact, they show that $c(X)=hc(X)=hL(X)$ for $X$ monotonically
normal.  Gartside's paper \cite{gartside} contains an extensive treatment of cardinal invariants of monotonically normal spaces.}

The other alternative available to us is that there is an infinite $A$ such that $\{K_\alpha:\alpha\in A\}$
is pairwise disjoint. This case is disposed of by essentially the same argument used in the proof of Proposition~\ref{prop1},
so we leave it to the reader.

Since either alternative leads to a contradiction, it must be the case that there is an $\alpha<\omega_1$
for which all $=_\alpha$--classes are of cardinality at most $2$ and this establishes the proposition
\end{proof}

Finally, to show that in fact the conclusion of Proposition~\ref{prop2} holds for {\em every} such model
$M$, we note that the definition of our equivalence relations $=_M$ doesn't depend on the fact that $M$
is an elementary submodel of $H(\chi)$.  In fact, we can carry out the same construction
given  {\em any} set of continuous real--valued functions defined on $X$.  Thus, we can view Proposition~\ref{prop2}
as stating that there is a countable set $\mathcal{X}$ of continuous real--valued functions defined on $X$
for which the associated equivalence classes all have cardinality at most 2  --- simply take $\mathcal{X}$
to be the set of all such functions in $M$.

Now let $M$ be an {\em arbitrary} countable elementary submodel of $H(\chi)$ containing $X$.  By elementarity,
\begin{equation}
M\models\text{``there is a countable set $\mathcal{X}$ as above''.}.
\end{equation}
Since $\mathcal{X}$ is countable, every member of $\mathcal{X}$ is also in $M$.  Thus, any $=_M$--equivalence
class is contained in a $=_{\mathcal{X}}$--equivalence class, and therefore each such $=_M$--class is of
cardinality at most 2.  Thus, the proof of Theorem~\ref{mainthm} is complete.
\end{proof}

In closing, we remark that the rest of Rudin's proof of Nikiel's Conjecture seems to be amenable to a similar treatment, and we
plan to examine this in a future paper.

\bibliographystyle{plain}

\end{document}